\documentstyle[12pt]{amsart}
\newtheorem{theorem}{Theorem}[section]
\newtheorem{corollary}[theorem]{Corollary}

\numberwithin{equation}{section}

\begin{document}
\title [Symplectic mean curvature flow]
 {Symplectic mean curvature flows in K\"ahler surfaces with positive holomorphic sectional curvatures}
\author{Jiayu Li, Liuqing Yang}

\address{Jiayu Li, School of Mathematical Sciences, University of Science and Technology of China Hefei 230026 \\ AMSS CAS Beijing 100190, P. R. China}
\email{lijia@@amss.ac.cn}

\address{
Liuqing Yang, Academy of Mathematics and Systems Sciences\\
Chinese Academy of Sciences\\ Beijing 100190, P. R. of China.}
\email{yangliuqing@@amss.ac.cn}


\date{}

\maketitle

\noindent\textbf{Abstract}: In this paper, we mainly study the
mean curvature flow in K\"ahler surfaces with positive holomorphic
sectional curvatures. We prove that if the ratio of the maximum
and the minimum of the holomorphic sectional curvatures is less
than $2$, then there exists a positive constant $\delta$ depending
on the ratio such that $\cos\alpha\geq\delta$ is preserved along
the flow.

\vspace{.2in}


\section{Introduction}

\allowdisplaybreaks

\vspace{.1in}

\noindent Mean curvature flows were studied by many authors, for
example Huisken (\cite{H1}, \cite{H2}), Ecker-Huisken (\cite{EH}),
Huisken-Sinestrari \cite{HS}, Carlo Ilmanen \cite{I}, Neves
\cite{N}, Smoczyk (\cite{Sm1}), Wang (\cite{W}), White \cite{Wa},
etc.

In this paper we mainly concentrated on the symplectic mean
curvature flows, which were studied by Chen -Tian \cite{CT},
Chen-Li \cite{CL1}, Chen-Li-Tian \cite{CLT}, Wang \cite{W}, Han-Li
\cite{HL}, \cite{HL3}, \cite{HL2}, Han-Sun \cite{HS}, Han-Li-Sun
\cite{HLS}, and Han-Li-Yang \cite{HLY}. The basic fact is that the
symplectic property is preserved by the mean curvature flow if the
ambient space $M$ is K\"ahler-Einstein, or if the ambient K\"ahler
surface evolves along the K\"ahler-Ricci flow \cite{HL2}.

Let $(M,J,\overline\omega,\bar{g})$ be a K\"ahler surface. For a
compact oriented real surface $\Sigma$ which is smoothly immersed
in $M$, the K\"ahler angle \cite{CW} $\alpha$ of $\Sigma$ in $M$
was defined by

$$\omega|_\Sigma=\cos\alpha d\mu_\Sigma$$ where $d\mu_\Sigma$ is
the area element of $\Sigma$ in the induced metric from
$\overline{g}$. We say that $\Sigma$ is a symplectic surface if
$\cos\alpha > 0$; $\Sigma$ is a holomorphic curve if $\cos\alpha
\equiv1$.

Given an immersed $F_0: \Sigma\rightarrow M$, we consider a
one-parameter family of smooth maps $F_t=F(\cdot, t):
\Sigma\rightarrow M$ with corresponding images
$\Sigma_t=F_t(\Sigma)$ immersed in $M$ and $F$ satisfies the mean
curvature flow equation:
\begin{equation}\label{flow}
   \begin{cases}
     \frac{\partial}{\partial t} F(x,t) =H(x,t) \\
      F(x,0)=F_{0}(x),
   \end{cases}
\end{equation}
where $H(x,t)$ is the mean curvature vector of $\Sigma_t$
at $F(x,t)$ in $M$.

Choose an orthonormal basis $\{e_1,e_2,e_3,e_4\}$ on ($M$,$\bar g$)
along $\Sigma_t$ such that $\{e_1,e_2\}$ is the basis of $\Sigma_t$
and the symplectic form $\omega_t$ takes the form
\begin{eqnarray}\label{omega}
  \omega_t=\cos\alpha u_1\wedge u_2+\cos\alpha u_3\wedge
  u_4+\sin\alpha u_1\wedge u_3-\sin\alpha u_2\wedge u_4,
\end{eqnarray}
where $\{u_1,u_2,u_3,u_4\}$ is the dual basis of
$\{e_1,e_2,e_3,e_4\}$. Then along the surface $\Sigma_t$ the complex
structure on $M$ takes the form (\cite{CL1})
\begin{eqnarray}\label{J}
J= \begin{pmatrix}
0 & \cos\alpha & \sin\alpha & 0\\
-\cos\alpha & 0 & 0 & -\sin\alpha\\
-\sin\alpha & 0 & 0 & \cos\alpha \\
0 & \sin\alpha & -\cos\alpha & 0
\end{pmatrix}.\end{eqnarray}

Recall the evolution equation of the K\"ahler angle along the mean
curvature flow deduced in \cite{HL2},
\begin{theorem}\label{angle}
  The evolution equation for $\cos\alpha$ along $\Sigma_t$ is
  \begin{eqnarray}\label{alpha}
  (\frac{\partial}{\partial t}-\Delta)\cos\alpha&=&|\overline\nabla
  J_{\Sigma_t}|^2\cos\alpha+\sin^2\alpha Ric(Je_1,e_2).
  \end{eqnarray}
\end{theorem}
Here
\begin{eqnarray}
  |\overline\nabla
  J_{\Sigma_t}|^2=|h^4_{1k}+h^3_{2k}|^2+|h^4_{2k}-h^3_{1k}|^2\geq\frac{1}{2}|H|^2.
\end{eqnarray}
 We want to see whether the symplectic property is preserved along the mean
curvature flow. In the case that $M$ is a K\"ahler-Einstein
surface, we have $Ric(Je_1,e_2)=\bar{\rho}\cos\alpha$, where
$\bar{\rho}$ is the scalar curvature of $M$, so the symplectic
property is preserved. If the ambient K\"ahler surface evolves
along the K\"ahler-Ricci flow, Han-Li \cite{HL2} derived the
evolution equation for $\cos\alpha$ and consequently they showed
that the symplectic property is also preserved. In this paper, we
find another condition to assure that along the flow, at each time
the surface is symplectic. Note that we don't require $M$ to be
Einstein. Denote the minimum and maximum of holomorphic sectional
curvatures of $M$ by $k_1$ and $k_2$. We state our main theorem as
follows:

\vspace{.1in} \noindent {\bf Main Theorem} {\it Suppose $M$ is a
K\"ahler surface with positive holomorphic sectional curvatures. Set
$\lambda=\frac{k_2}{k_1}$. If the flow satisfies either

  I. $1\leq\lambda<\frac{11}{7}$ and
  $\cos\alpha(\cdot,0)\geq\delta>\frac{53(\lambda-1)}{\sqrt{(53\lambda-53)^2+(48-24\lambda)^2}}$,

  or

  II. $\frac{11}{7}\leq\lambda<2$ and
$\cos\alpha(\cdot,0)\geq\delta>\frac{8\lambda-5}{\sqrt{(8\lambda-5)^2+(12-6\lambda)^2}}$,

\noindent then along the flow
  \begin{eqnarray}\label{my alpha}
    (\frac{\partial}{\partial t}-\Delta)\cos\alpha\geq|\overline\nabla
    J_{\Sigma_t}|^2\cos\alpha+C\sin^2\alpha,
  \end{eqnarray}
  where $C$ is a positive constant depending only on $k_1$,
  $k_2$ and $\delta$. As a corollary, $\min_{\Sigma_t}\cos\alpha$ is increasing with respect to $t$. In particular, at each time $t$, $\Sigma_t$ is
  symplectic. Therefore, we call this flow the symplectic mean curvature flow.}

Since we obtain (\ref{my alpha}), many theorems in ``symplectic mean
curvature flows in K\"ahler-Einstein surfaces" still hold in our
case. For example,

 Arguing as in \cite{CW} by strong maximum principle, we have
\begin{corollary}
I. Suppose $M$ is a K\"ahler surface with positive holomorphic
sectional curvatures and $1\leq\lambda<\frac{11}{7}$, then every
symplectic minimal surface satisfying
\[\cos\alpha>\frac{53(\lambda-1)}{\sqrt{(53\lambda-53)^2+(48-24\lambda)^2}}\]
in $M$ is a holomorphic curve.

II. Suppose $M$ is a K\"ahler surface with positive holomorphic
sectional curvatures and $\frac{11}{7}\leq\lambda<2$, then every
symplectic minimal surface satisfying
\[\cos\alpha>\frac{8\lambda-5}{\sqrt{(8\lambda-5)^2+(12-6\lambda)^2}}\]
in $M$ is a holomorphic curve.
\end{corollary}
Arguing exactly in the same way as in \cite{CL1} or \cite{W}, we have
\begin{theorem}\label{singularity}
Under the same condition of the Main Theorem, the symplectic mean
curvature flow has no type I
  singularity at any $T>0$.
\end{theorem}

Acknowledgment: The research was supported by NSFC 11071236,
11131007 and 10421101.

\section{Curvature Tensor, Sectional Curvature and Holomorphic Sectional Curvature}
Denote the curvature tensor of $M$ by $K$. Set $K(X)=K(X,JX,X,JX)$
and $K(X,Y)=K(X,Y,X,Y)$, where $X,Y$ are arbitrary vector fields on
$M$. It is known that (c.f. \cite{BGG}, \cite{V}) we can express the sectional curvatures by holomorphic
sectional curvatures.
\begin{theorem}
The sectional curvatures of $M$ can be determined by
the holomorphic sectional curvatures by
  \begin{eqnarray}\label{sectional}
      K(X,Y)&=&\frac{1}{32}[3K(X+JY)+3K(X-JY)-K(X+Y)-K(X-Y)\nonumber\\
      &&-4K(X)-4K(Y)].
  \end{eqnarray}
\end{theorem}

Using (\ref{sectional}), it is easy to check that,
\begin{theorem} For any vector fields $X$, $Y$ and $Z$ on $M$,
\begin{eqnarray}\label{xyxz}
  &&K(X,Y,X,Z)\nonumber\\
  &=&\frac{1}{2}[K(Y+Z,X)-K(X,Y)-K(X,Z)]\nonumber\\
  &=&\frac{1}{64}[3K(Y+Z+JX)+3K(Y+Z-JX)-K(Y+Z+X)-K(Y+Z-X)\nonumber\\
  &&-3K(Y+JX)-3K(Y-JX)-3K(Z+JX)-3K(Z-JX)\nonumber\\
  &&-4K(Y+Z)+K(Y+X)+K(Y-X)+K(Z+X)+K(Z-X)\nonumber\\
  &&+4K(X)+4K(Y)+4K(Z)].
  \end{eqnarray}
\end{theorem}

Denote the minimum and the maximum of sectional curvatures by
$K_{min}$ and $K_{max}$ respectively, we have the following
estimates.
\begin{theorem}\label{K}$K_{min}$ and $K_{max}$ satisfy
\begin{eqnarray}
  K_{max}\leq\frac{3}{2}k_2-\frac{1}{2}k_1
  \end{eqnarray}
  and
  \begin{eqnarray}
    K_{min}\geq\frac{3}{4}k_1-\frac{1}{2}k_2
  \end{eqnarray}
\end{theorem}
{\it Proof.} Given any point $p\in M$ and any two unit orthogonal
vectors $X$ and $Y$ at $p$, we can find two vectors $Z$ and $W$ such
that $\{X,Y,Z,W\}$ form an orthonormal basis of $T_pM$. Suppose
$JX=yY+zZ+wW$, then
\begin{eqnarray}{\label{X+JY}}
  \langle X+JY,X+JY\rangle=2-2y,
\end{eqnarray}
and
\begin{eqnarray}{\label{X-JY}}
  \langle X-JY,X-JY\rangle=2+2y.
\end{eqnarray}

 Assume the K\"ahler form is anti-self-dual, it was shown in
\cite{HLY} that, $y^2+z^2+w^2=1$ and $J$ has the form
\begin{eqnarray}\label{J1}
J=\left (\begin{array}{clcr} 0 &y &z &w \\
-y &0 &w &-z\\
-z &-w &0 &y\\
-w &z &-y &0 \end{array}\right).
\end{eqnarray}
Combining (\ref{sectional}) with (\ref{X+JY}) and (\ref{X-JY}), we
get
\begin{eqnarray*}
  K(X,Y)&\leq&\frac{1}{32}[3(2-2y)^2k_2+3(2+2y)^2k_2-2^2k_1-2^2k_1-4k_1-4k_1]\\
  &=&\frac{1}{4}[(3+3y^2)k_2-2k_1]\\
  &\leq&\frac{1}{4}(6k_2-2k_1)\\&=&\frac{3}{2}k_2-\frac{1}{2}k_1,
\end{eqnarray*}
and similarly
\begin{eqnarray*}
  K(X,Y)&\geq&\frac{1}{4}[(3+3y^2)k_1-2k_2]\\
  &\geq&\frac{1}{4}(3k_1-2k_2)\\&=&\frac{3}{4}k_1-\frac{1}{2}k_2.
\end{eqnarray*}
This proves the theorem.   \hfill Q.E.D.

\vspace{.2in}

\section{Proof of the Main Theorem}
In this section, we will prove the Main Theorem of this paper.

\vspace{.2in}

 {\it Proof of the Main Theorem.} In order to prove
this theorem, we need to estimate $Ric(Je_1,e_2)$. Using two
different methods, we get two available estimates. We now deduce
the first one.
\begin{eqnarray}\label{Ric}
  Ric(Je_1,e_2)&=&K(Je_1,e_1,e_2,e_1)+K(Je_1,e_3,e_2,e_3)+K(Je_1,e_4,e_2,e_4)\nonumber\\
  &=&K(\cos\alpha e_2+\sin\alpha e_3,e_1,e_2,e_1)+K(\cos\alpha e_2+\sin\alpha
  e_3,e_3,e_2,e_3)\nonumber\\
  &&+K(\cos\alpha e_2+\sin\alpha e_3,e_4,e_2,e_4)\nonumber\\
  &=&cos\alpha R_{22}+\sin\alpha(K_{3121}+K_{3424}),
\end{eqnarray}
where
\begin{eqnarray}\label{22}
  R_{22}=K_{2121}+K_{2323}+K_{2424}.
\end{eqnarray}
By (\ref{sectional}), we have
\begin{eqnarray*}
  K_{2121}&=&\frac{1}{32}[3K(e_1+Je_2)+3K(e_1-Je_2)-K(e_1+e_2)-K(e_1-e_2)\\
  &&-4K(e_1)-4K(e_2)].
\end{eqnarray*}
By our choice of the complex structure (\ref{J}), we get
\begin{eqnarray*}
  \langle e_1+Je_2,e_1+Je_2\rangle=2-2\cos\alpha,
\end{eqnarray*}
and
\begin{eqnarray*}
  \langle e_1-Je_2,e_1-Je_2\rangle=2+2\cos\alpha.
\end{eqnarray*}
Hence $K_{2121}$ can be estimated by $k_1$ and $k_2$,
\begin{eqnarray}\label{2121}
  K_{2121}&\geq&\frac{1}{32}[3(2-2\cos\alpha)^2k_1+3(2+2\cos\alpha)^2k_1-2^2k_2-2^2k_2-4k_2-4k_2]\nonumber\\
  &=&\frac{1}{4}[(3+3\cos^2\alpha)k_1-2k_2].
\end{eqnarray}
Similarly, we get
\begin{eqnarray}\label{2323}
  K_{2323}\geq\frac{1}{4}(3k_1-2k_2),
\end{eqnarray}
and
\begin{eqnarray}\label{2424}
  K_{2424}\geq\frac{1}{4}[(3+3\sin^2\alpha)k_1-2k_2].
\end{eqnarray}
Putting (\ref{2121}), (\ref{2323}) and (\ref{2424}) into (\ref{22}),
we obtain that
\begin{eqnarray}\label{R22}
  R_{22}\geq3k_1-\frac{3}{2}k_2.
\end{eqnarray}
Using (\ref{xyxz}) and (\ref{J}), we can also estimate $K_{3121}$
and $K_{3424}$. We have
\begin{eqnarray}\label{3121}
K_{3121}\geq\frac{1}{32}[(53+48\sin\alpha\cos\alpha)k_1-53k_2],
\end{eqnarray}
and
\begin{eqnarray}\label{3424}
K_{3424}\geq\frac{1}{32}[(53-48\sin\alpha\cos\alpha)k_1-53k_2].
\end{eqnarray}
Adding (\ref{3121}) and (\ref{3424}) yields
\begin{eqnarray}
  K_{3121}+K_{3424}\geq-\frac{53}{16}(k_2-k_1).
\end{eqnarray}
By a similar computation in the opposite direction, we get
\begin{eqnarray}\label{31213424}
  |K_{3121}+K_{3424}|\leq\frac{53}{16}(k_2-k_1).
\end{eqnarray}
 Therefore by (\ref{Ric}), (\ref{R22}), (\ref{31213424}) and short time existence of the mean curvature flow, we have
\begin{eqnarray}\label{positive}
  Ric(Je_1,e_2)&\geq&\cos\alpha(3k_1-\frac{3}{2}k_2)-\sqrt{1-\cos^2\alpha}\frac{53}{16}(k_2-k_1)\nonumber\\
  &=&(3\cos\alpha+\frac{53}{16}\sqrt{1-\cos^2\alpha})k_1-(\frac{3}{2}\cos\alpha+\frac{53}{16}\sqrt{1-\cos^2\alpha})k_2\nonumber.\\
\end{eqnarray}
If $1\leq\lambda<2$ and
$\cos\alpha>\frac{53(\lambda-1)}{\sqrt{(53\lambda-53)^2+(48-24\lambda)^2}}$,
then the RHS of (\ref{positive}) is positive.

Another estimate follows directly from Theorem \ref{K} and
Berger¡¯s inequality (c.f. \cite{G}) that
\begin{eqnarray}|K_{3121}+K_{3424}|\leq|K_{3121}|+|K_{3424}|\leq
K_{max}-K_{min}\leq 2k_2-\frac{5}{4}k_1.\end{eqnarray} Putting the
above estimate into $(\ref{Ric})$ yields
\begin{eqnarray}\label{positive2}
  Ric(Je_1,e_2)&\geq&\cos\alpha(3k_1-\frac{3}{2}k_2)-\sqrt{1-\cos^2\alpha}(2k_2-\frac{5}{4}k_1)\nonumber\\
  &=&(3\cos\alpha+\frac{5}{4}\sqrt{1-\cos^2\alpha})k_1-(\frac{3}{2}\cos\alpha+\frac{5}{4}\sqrt{1-\cos^2\alpha})k_2\nonumber.\\
\end{eqnarray}
It follows that if $1\leq\lambda<2$ and
$\cos\alpha>\frac{8\lambda-5}{\sqrt{(8\lambda-5)^2+(12-6\lambda)^2}}$,
then the RHS of (\ref{positive}) is positive. Note that
\begin{eqnarray*}\frac{53(\lambda-1)}{\sqrt{(53\lambda-53)^2+(48-24\lambda)^2}}\leq\frac{8\lambda-5}{\sqrt{(8\lambda-5)^2+(12-6\lambda)^2}}
\end{eqnarray*}
for $1\leq\lambda<\frac{11}{7}$, and
\begin{eqnarray*}\frac{53(\lambda-1)}{\sqrt{(53\lambda-53)^2+(48-24\lambda)^2}}\geq\frac{8\lambda-5}{\sqrt{(8\lambda-5)^2+(12-6\lambda)^2}}
\end{eqnarray*}
for $\frac{11}{7}\leq\lambda<2$, we get the conclusion.   \hfill
Q.E.D.

\end{document}